\documentclass[article,12pt]{article}
\usepackage{graphicx}
\usepackage{amsmath}
\usepackage{mathrsfs,amsfonts}
\usepackage{amsfonts}
\usepackage{amsmath,amssymb,amsfonts}
\usepackage{amssymb}
\usepackage{color}
\usepackage{amsmath}  
\usepackage{amsfonts} 
\usepackage{graphicx} 
\usepackage{braket}
\usepackage{enumerate}
\usepackage{indentfirst}
\usepackage{longtable}
\usepackage{setspace}
\usepackage{authblk}

\allowdisplaybreaks[4]

\def\al{{\alpha}}\def\be{{\beta}}\def\de{{\delta}}
\def\ep{{\epsilon}}\def\ka{{\kappa}}
\def\si{{\sigma}}
\def\ze{{\zeta}}

\def\th{{\theta}}

\def\<{\left<}\def\>{\right>}\def\({\left(}\def\){\right)}

\newfam\msbmfam\font\tenmsbm=msbm10\textfont
\msbmfam=\tenmsbm\font\sevenmsbm=msbm7

\scriptfont\msbmfam=\sevenmsbm\def\bb#1{{\fam\msbmfam #1}}

\def\BB{\bb B}
\def\EE{\bb E}\def\FF{\bb F}
\def\PP{\bb P}
\def\RR{\bb R}

\usepackage{amsmath,amsfonts,amsthm,amssymb}    
\usepackage{gensymb, mathrsfs}                  

\usepackage{accents}
\DeclareMathSymbol{\widehatsym}{\mathord}{largesymbols}{"62}

\def\*#1{\mathbf{#1}}

\renewcommand{\tilde}{\widetilde}
\renewcommand{\phi}{\varphi}

\renewcommand{\qed}{\hfill \ensuremath{\Box}}

\makeatletter \renewcommand\d[1]{\ensuremath{%
  \;\mathrm{d}#1\@ifnextchar\d{\!}{}}}
\makeatother

\theoremstyle{plain}

\theoremstyle{definition}


\newcommand{\beq}{\begin{equation}}
\newcommand{\eeq}{\end{equation}}

\definecolor{c}{rgb}{0.9,0.3,0.1}
\definecolor{b}{rgb}{0.1,0.3,0.9}

\newtheorem{remark}{Remark}[section]
\newtheorem{lemma}{Lemma}[section]
\newtheorem{theorem}{Theorem}[section]
\newtheorem{definition}{Definition}[section]
\newtheorem{hypothesis}{Hypothesis}[section]

\renewcommand{\theequation}{\arabic{section}.\arabic{equation}}

\def\al{{\alpha}}\def\be{{\beta}}\def\de{{\delta}}
\def\ep{{\epsilon}}

\def\si{{\sigma}}\def\th{{\theta}}\def\ze{{\zeta}}

\def\<{\left<}\def\>{\right>}\def\({\left(}\def\){\right)}

\newfam\msbmfam\font\tenmsbm=msbm10\textfont
\msbmfam=\tenmsbm\font\sevenmsbm=msbm7
\scriptfont\msbmfam=\sevenmsbm\def\bb#1{{\fam\msbmfam #1}}
\def\BB{\bb B}
\def\EE{\bb E}\def\PP{\bb P}
\def\RR{\bb R}

\numberwithin{equation}{section}

\begin{document}

\title{\large \bf Stochastic maximum principle for weighted mean-field system}

\author{Yanyan Tang\thanks{ Department of Mathematics, Southern University of Science and Technology, Shenzhen, Guangdong, 518055, China ({\tt tangxiaotang9600@163.com})} and 
Jie Xiong\thanks{(Corresponding author. Shenzhen Key Laboratory of Safety and Security for Next Generation of Industrial Internet and  Southern University of Science and Technology, Shenzhen, Guangdong, 518055, China {\tt xiongj@sustech.edu.cn}).
This author is supported by ZDSYS20210623092007023 and
 NSFC Grants 61873325 and 11831010.}
}
\date{}
\maketitle
 \bigskip

\noindent \textbf{Abstract.}
We study the optimal control problem for a weighted mean-field system. A new feature of the control problem is that the coefficients   depend on the state  process as well as  its weighted measure and the control variable. By applying variational technique, we establish a stochastic maximum principle. As an application, we investigate the optimal premium policy of an insurance firm for asset–liability management problem. 
\bigskip

\noindent \textbf{Keyword.}
 McKean-Vlasov equation,  stochastic maximum principle.\\
\noindent \textbf{AMS subject classifications.}
60H10, 93E20,  93E03, 60H30

\section{Introduction}\label{intro} \label{sec1}
\setcounter{equation}{0}
\renewcommand{\theequation}{\thesection.\arabic{equation}}

Optimal control problems have been studied extensively since the pioneering work of Pontryagin \cite{P-1987}, where a maximum principle  is obtained by  using spike variation. Kushner (\cite{K-1965}, \cite{K-1972})  investigated the applicability of the maximum principle to the design of controllers for stochastic system, i.e. stochastic maximum principle (SMP) for optimal control. He  showed that for system with additive white noise a maximum principle of the form developed by Pontryagin is valid  when diffusion coeffcient does not depend on the control variable. When the diffusion coefficient contains a control variable, Bensoussan (\cite{B-1981}, \cite{B-1984}) studied such a case. The maximum principle he obtained are local conditions, and his method depends heavily on the control domain being convex. Peng \cite{Peng-1990} broke through this difficulty in 1990. In that paper, first and second order variational inequalities are introduced, when the control domain need not to be convex, and the diffusion coefficient contains the control variable.

In recent years, stochastic optimal control problems for the mean-field stochastic differential equations (SDEs) have attracted an increasing attention. The mean-field SDEs can trace their roots to the McKean–Vlasov model, which was first introduced by Kac \cite{K-1956} and McKean \cite{M-1966}  to study physical systems with a large number of interacting particles. Lasry and Lions \cite{2009+} extended the applications of the mean-field models to economics and finance.  It is not until Buckdahn et al \cite{B-L-P-2009}  and Buckdahn et al \cite{2009+} established the theory of the mean-field BSDEs that the SMP for the optimal control system of mean-field type has become a popular topic. For example, Li \cite{L-2012} studied SMP for mean-field controls, with the domain of the control convex; Buckdahn et al \cite {B-L-M-2016-} obtained the related SMP for a class of general stochastic control problem with McKean-Vlasov dynamics, in which the coefficients depend nonlinearly on both the state process as well as its distribution; Acciaio et al. \cite{A-2018} studied mean-field stochastic control problems where the cost functional and the state dynamics depend on the joint distribution of the controlled state and the control process. Strongly inspired by Buckdahn et al \cite{B-L-M-2016-}, Buckdahn et al \cite{B-C-L-2021} investigated a generalized mean-field SMP for the optimal control problem where the coefficients not only depend on the distribution of $(X,v)$ but also with partial information. Interested readers may refer to Shen and Siu \cite{S-S-2013},  Djehiche et al \cite{D-T-2015}, Guo and Xiong \cite{G-X-2018}, Zhang et al \cite{Z-S-X-2018}, Lakhdari et al \cite{L-M-H-2021}, Zhang \cite{Z-2021} and Wang and Wu \cite{W-W-2022} for various versions of the SMPs for the mean-field models.

As we mentioned above, one of the motivations of the study of mean-field SDE is from mathematical finance. The state process is usually the valuation of a certain asset, such as the stock price, by a typical investor which is affected by  this investor as well as the others who are interested in this asset, in a collective manner. So far, this collective interaction is represented by the average measure of the whole system. By a propagation of chaos limit, it is described by the distribution of the state process. However, in the real-world situations, the opinion of the investors about the asset should be weighted accord to their wealth levels. This motivates our study of optimal control with weighted mean-field interaction.

In this paper, we are intersted in the SMP for a stochastic control problem with McKean-Vlasov dynamics: 
\begin{equation}\label{eq0613a-}
\left\{\begin{array}{ccl}
dX_t&=&b(X_t,\mu^{X,A}_t,u_t)dt+\si(X_t,\mu^{X,A}_t,u_t)dW_t\\
dA_t&=&A_t\(\al(X_t,\mu^{X,A}_t,u_t)dt+\be(X_t,\mu^{X,A}_t,u_t)dW_t\)  \ \ \ t\in[0,T] \\
X_0&=&x,\ A_0=a,   
\end{array}\right.\end{equation}
where $W_t$ is an $m$-dimension Browanion motion defined on a complete probability space $(\Omega
,\mathscr{F},\PP)$. In the above, $(b,\si):\ \RR^d\times\mathcal{M}_F(\RR^d)\times U\rightarrow \RR^d\times \RR^{d\times m}$ and
$(\al,\be):\ \RR^d\times\mathcal{M}_F(\RR^d)\times U\rightarrow \RR\times \RR^{1\times m}$ are continous mappings, $u_\cdot\in L_{\FF}^2([0,T];U)$ is a control process, $U\subset\RR^k$ is a convex domain for the control, and $\mu^{X,A}_t\in\mathcal{M}_F(\RR^d)$ is the weighted measure given by
\begin{equation}\label{weighted measure}
\<\mu^{X,A}_t,f\>=\EE\(A_t f(X_t)\),\qquad\forall f\in C_b( \RR^d).
\end{equation}
The equation of the form (\ref{eq0613a-}) without the control has been studied by Kurtz and Xiong \cite{K-X-1999} as the limit of a system of interacting weighted particles. 
The main difficulty in obtaining the uniqueness of the solution to this equation as well as the interacting system related to it is the non-Lipschitz property of the coefficients such as $(a,x,\mu)\mapsto a\al(x,\mu,u)$ because $a$ is unbounded. A stopping time technique was used in \cite{K-X-1999}. However, for stochastic control problem, the adjoint equation will be a backward stochastic differential equation (BSDE), and hence, the stopping argument will not be convenient. In this article, we overcome this difficulty by estimating $X_t$ and $A_t$ using suitably chosen different norms.

The rest of this paper is organized as follows: In the next section, we  formulate the problem and present the stochastic maximum principle for  our stochastic control problem. In Section 3, as an application, we give  an example about the optimal premium policy for  an insurance firm for asset–liability management problem.
In Sections 4 and 5, we give detailed proofs about the existence and  uniqueness for the solution  of  the  weighted mean-field system and about the 
stochastic maximum principle, respectively.

\section{Stochastic Maximum Principle}

We recall that we cosider the stochastic control problem with the state equation (\ref{eq0613a-}). Let ${{\mathcal{U}}}[0,T]$ be the collection of all $U$-valued ${\mathcal F}_t$-adapted processes $u_t$ satisfying
\[\EE \int_0^T|u_t|^2dt<\infty.\]
${{\mathcal{U}}}[0,T]$ is called the set of admissible controls. In order not to over complicate the already notational heavy presentation
of this paper, in what follows we shall assume all processes are 1-dimensional (i.e., $d = k = m = 1$). We should note that the higher dimensional cases can be argued along the same lines without substantial difficulties, except for even heavier notations.

We will take the cost function $J(u)$ as 
\begin{equation}\label{eq0613c1}
J(u)=\EE\(\int^T_0f(X_t,A_t,\mu^{X,A}_t,u_t)dt+\Phi(X_T,A_T)\),
\end{equation}
where $f:\RR\times\RR\times\mathcal{M}_F(\RR)\times{U}\rightarrow\RR$, $\Phi:\RR\times\RR\rightarrow\RR.$ Here we have dropped 
the superscript $u$ of $X,A$ for simplicity of notation.

Our goal is to find a control to minimize the cost functional over $\mathcal{U}$, namely, an admissible control $u\in\mathcal{U}$ is said to the optimal if 
\begin{equation}\label{eq0613c2}
J(u)=\min_{v\in\mathcal{U}}J(v).
\end{equation}

In this section, we consider the necessary condition for the optimal control of  the problem (\ref{eq0613c2}), also know as stochastic maximum principle,  with  the convexity assumption on the control set ${U}$.

 For weighted mean-field type SDEs and BSDEs, we have to introduce some notations. Let $(\Omega',\mathscr{F}',\PP')$ be
  a copy of the probability space $(\Omega,\mathscr{F},\PP)$ and  $\EE'(\cdot)=\int_{\Omega'}(\cdot)d\PP'$, $\xi'$ is a random variable defined on $(\Omega',\mathscr{F}',\PP')$. For $\mu_1,\mu_2\in\mathcal{M}_F(\RR),$ the Wasserstein metric is defined by 
\[\rho(\mu_1,\mu_2)=\sup\{|\<\mu_1,f\>-\<\mu_2,f\>|: f\in\BB_1\},\]
where $\BB_1=\{f:|f(x)-f(y)|\leq|x-y|,|f(x)|\leq1,\forall x,y\in\RR\}$, $\<\mu,f\>$ stands for the integral of the function $f$ with respect to the measure $\mu$.

\begin{definition}
Suppose $f:\mathcal{M}_F(\RR)\to\RR$. We say that $f\in C^1(\mathcal{M}_F(\RR))$ if there exists $h(\mu;\cdot)\in C_b(\RR)$ such that
\[f(\mu+\ep\nu)-f(\mu)=\<\nu,h(\mu;\cdot)\>\ep+o(\ep).\]
We denote $h(\mu;x)=f_\mu(\mu;x)$. 
\end{definition}

Next, we make the following standard assumptions.
\begin{hypothesis}\label{hyp}
The conefficients $b,\sigma,\al, \be, f,\Phi$ are measurable in all variables. Furthermore 
 \begin{enumerate}[(1)]
\item $b,\sigma,f\in C^1(\RR\times \mathcal{M}_F(\RR)\times U)$ with bounded partial derivatives;
\item\label{2} $\al,\be\in C^1(\RR\times \mathcal{M}_F(\RR)\times U)$ are bounded and  with bounded partial derivatives; 
\item $\phi_\mu(x,\mu,u,x')$ is differentable in $x'$ with bounded derivative for $\phi=b,\sigma,\al,\be$. We denote the partial derivative of  $\phi_\mu(x,\mu,u;x')$ with respect to $x'$ by $\phi_{\mu,1}(x,\mu,u;x')$;
\item $\Phi(x,a)$ has bounded partial derivatives in $x,a$.
\end{enumerate}
\end{hypothesis}

\begin{theorem}\label{thm0730a}
 For  $u\in\mathcal{U}$ being fixed, the  weighted mean-field SDE (\ref{eq0613a-}) has a unique solution.
\end{theorem}

The main result of this paper is to prove the following SMP. More precisely, we define Hamiltonnian
\begin{eqnarray*}
H(X,A,\mu,u,P,Q,p,q)&:=&A\(p\al(X,\mu,u)+q\be(X,\mu,u)\)\\
&&+Pb(X,\mu,u)+Q\si(X,\mu,u)-f(X,A,\mu,u).
\end{eqnarray*}

The  adjoint processes $(p_t,q_t),(P_t,Q_t)$ are governed by the following equations 
\begin{equation}\label{eq0613adjointa1}
\left \{
\begin{split}
  \begin{aligned}
   dp_t=&-\big\{p_t\al(\th_t)+q_t\be(\th_t)-f_a(\ka_t)\big\}dt&\\
		    &-\mathbb{E}' \big\{\(P_t'b_{\mu}+Q_t'\si_{\mu}+A_t'(p_t'\al_{\mu}+q_t'\be_{\mu})\)(\th_t';X_t)-f_\mu(\ka_t';X_t)\big\}dt+q_t dW_t,&\\
   p_T=&-\Phi_a(X_T,A_T),&
      \end{aligned}
  \end{split}
  \right.
\end{equation}
and 
\begin{equation}\label{eq0613adjointa2}
\left \{
\begin{split}
  \begin{aligned}
   dP_t=&-\big\{P_tb_x(\th_t)+Q_t\si_x(\th_t)+A_t(p_t\al_x(\th_t)+q_t\be_x(\th_t))-f_x(\ka_t)\big\}dt&\\
    &-\mathbb{E}' \big\{A_t\(\big(P_t'b_{\mu1}+Q_t'\si_{\mu1}+A_t'(p_t'\al_{\mu1}+q_t'\be_{\mu1})\)(\th_t';X_t)-f_{\mu1}(\ka_t';X_t)\big)\big\}dt&\\
&+Q_t dW_t,&\\
   P_T=&-\Phi_x(X_T,A_T),&
      \end{aligned}
  \end{split}
  \right.
\end{equation}
where $\ka_t=(X_t,A_t,\mu_t,u_t), \th_t=(X_t,\mu_t,u_t)$, $\phi_x$ stands for the partial derivative of $\phi$ with respect to $x$. 
The same convention also apply to variables $a$ and $u$. Note that the subscript $t$, such as $X_t$, stands for the time parameter of a stochastic process.

We are now ready to state the main theorem of the paper.
\begin{theorem}\label{smp}
(Stochastic Maximum Principle) Suppose that Hypothesis \ref{hyp} holds. Let $(X_t,u_t)$ be an optional solution of control problem (\ref{eq0613c2}). Then there are two pairs of $\mathscr{F}_t$-adapted processes $(p_t,q_t)$ and $(P_t,Q_t)$ there satisfy the equations (\ref{eq0613adjointa1}) and (\ref{eq0613adjointa2}) such that 
\[H_u(X_t,A_t,\mu_t,u_t,P_t,Q_t,p_t,q_t)=0.\]
\end{theorem}

\section{An asset liability management problem}

As an application of Theorem \ref{smp}, in this section, we consider the optimal premium policy of an insurance firm for asset 
liability management problem. 

Let $u$ be the premium strategy of the firm; $X_t$ be the cash flow;  $l_t\equiv l_t^{X,A,u}$ be the liability process. Suppose that $ a_t,  b_t, c_t $ are deterministic and uniformly bounded on $[0, T ]$. Suppose that the liability process $l_t$ satisfies
\[dl_t=-\big(a_t\EE(A_tX_t)+b_tu_t\big)dt-c_tdW_t,\]
where $c_t$ denotes the volatility rate, and $A_t$ is the weighted process satisfying
\begin{equation}\label{weighted process }
dA_t=A_t\al_tdt+A_t\be_tdW_t
\end{equation}
with initial value $a_0$  and $\al_t,\be_t $ are deterministic and bounded on $[0,T]$. Suppose that the firm has no liability at time 0, and only invests in a money account due to certain market regulations. Accordingly, the insurance firm only invests in a money account with compounded interest rate $ \rho_t$, and hence its cash-balance process $X_t$ is 
\[X_t=e^{\int_0^t\rho_sds}\(X_0-\int_0^te^{-\int_0^s\rho_rdr}dl_t\),\ \ \ X_0=x_0\]
where $x_0\geq0$ represents the initial reserve. According to It$\hat{o}$'s formula, we see that the cash flow statisfies a mean-field SDE
\begin{equation}\label{eq0441 }
\left\{\begin{array}{ccl}
dX_t&=&\big(\rho_tX_t+a_t\EE(A_tX_t)+b_tu_t\big)dt+c_tdW_t,\\
X_0&=&x_0.
\end{array}\right.
\end{equation}
Let $U\subset\RR$ be convex. For any $u\in\mathcal{U}$, equation (\ref{eq0441 }) admits a unique adapted cash flow $X_t$. The insurance firm hopes
to drive its cash-balance process  to $c_0$ in average at the terminal time $T$ to meet some regulatory requirement, i.e. $\EE X_T=c_0$. In order to reconcile the contradiction between the liquidity  and profitability of insurance,  we introduce a performance functional of the firm
\begin{equation}\label{example-1-1}
J(u)=\frac{1}{2}\EE\big[\int_0^{T} (LX_t^2+MA_t^2+Nu_t^2)dt+ R( X_T-\EE(X_T))^2\big],
\end{equation}
with terminal constraint $ \EE(X_T)=c_0$, for some given $c_0\in\RR$, where $L, M, N, R$  are positive constants. 

Then the asset liability management problem is to find an  admissible premium policy $u\in\mathcal{U} $ such that  $J(u)=\inf_{v\in\mathcal{U}} J[v].$ Applying  Lagrange multiplier technique, we define a cost functional for any  $\lambda\in\RR$,
\begin{equation}\label{example-2}
J_\lambda(u)=\frac{1}{2}\EE\left[\int_0^{T} (LX_t^2+MA_t^2+Nu_t^2)dt+ R( X_T-c_0)^2+2\lambda(X_T-c_0)\right].
\end{equation}

For  every  $\lambda \in\RR$, we  find $u^\lambda$ to minimize $J_\lambda(u)$. In this situation, by Throrem \ref{smp} we write the Hamiltonian function 
\begin{eqnarray*}
H(X,A,\mu,u,P,Q,p,q)&=&A(\al_t p+\be_t q)+P\left(\rho_tX+a_t\int x\mu (dx)+b_tu\right)+c_tQ\\
&&-\frac{1}{2}\big(LX^2+MA^2+Nu^2\big),
\end{eqnarray*}
 and the adjoint equations
\begin{equation}\label{example-1}
\left \{
\begin{split}
  \begin{aligned}
   dp_t=&-\big(p_t\al_t+q_t\be_t-MA_t+a_tX_t\mathbb{E}(P_t) \big)dt+q_t dW_t,&\\
   p_T=&0,&
      \end{aligned}
  \end{split}
  \right.
\end{equation}
and 
\begin{equation}\label{e}
\left \{
\begin{split}
  \begin{aligned}
   dP_t=&-\big(\rho_tP_t-LX_t+a_tA_t\mathbb{E}(P_t)\big)dt+Q_t dW_t,&\\
   P_T=&-R(X_T-c_0)-\lambda.&
      \end{aligned}
  \end{split}
  \right.
\end{equation}
From Theorem \ref{smp}, we have 
\begin{equation}\label{control form-1}
Nu_t=b_tP_t.
\end{equation}
Hence, the related feedback control system takes the form
\begin{equation}\label{FBSDE-1}
\left\{\begin{array}{ccl}
dX_t&=&\big(\rho_tX_t+a_t\EE(A_tX_t)+N^{-1}b_t^2P_t\big)dt+c_tdW_t\\
dA_t&=&A_t\(\al_t dt+\be_t dW_t\)\\
dP_t&=&-\big(\rho_tP_t-LX_t+a_tA_t\mathbb{E}(P_t)\big)dt+Q_t dW_t,\\
X_0&=&x_0,\ A_0=a,\ P_T=-R(X_T-c_0)-\lambda, \ \ \ t\in[0,T].\\    
\end{array}\right.
\end{equation}
Here we have to deal with a fully coupled mean-field forward–backward SDE. In order to solve this equation we set 
\begin{equation}\label{P_t}
P_t=\phi_t X_t+\varphi_t\EE X_t+\chi_t\EE(A_tX_t)+\psi_t,
\end{equation}
where $\phi$ is  a  deterministic differentiable  function with terminal $-R$, $\varphi_t,\chi_t$, and $\psi_t$ are stochastic processes. We write 
\begin{equation}\label{eq0613 form}
\left\{\begin{array}{ccl}
d\varphi_t&=&\Gamma_t dt+\Lambda_tdW_t,\\
d\chi_t&=&\Pi_tdt+\Theta_tdW_t.\\  
d\psi_t&=&\Upsilon_tdt+\Sigma_tdW_t.\\ 
\varphi_T&=&\chi_T=0,\ \ \ \psi_T=Rc_0-\lambda.
\end{array}\right.
\end{equation}
Differentiating on $P_t$ in (\ref{P_t}) and comparing with (\ref{e}), we get
\begin{equation}\label{solve-2}
\left \{
\begin{split}
  \begin{aligned}
&-\rho_tP_t+LX_t-a_tA_t\EE (P_t)\\
&=\phi_t\rho_tX_t+\phi_ta_t\EE(A_tX_t)+\phi_tb_tu_t\\
&+\varphi_t\rho_t\EE(X_t)+\varphi_ta_t\EE(A_tX_t)+\varphi_tb_t\EE(u_t)\\
&+\chi_t(\rho_t+\al_t +a_t\EE (A_t))\EE(A_tX_t)+\chi_tb_t\EE(A_tu_t)+\be_t c_t\chi_t\EE(A_t)\\
&+\Upsilon_t+X_t\frac{d}{dt}\phi_t+\Gamma_t\EE (X_t) +\Pi_t\EE(A_tX_t)\\
&Q_t=\phi_tc_t+\Sigma_t+\Lambda_t\EE(X_t)+\Theta_t\EE(A_tX_t).
\end{aligned}
  \end{split}
  \right.
\end{equation}
From (\ref{control form-1}) and by comparing the coefficients of $X_t,\EE[X_t],\EE[A_tX_t],$   in the first equation of (\ref{solve-2}), we obtain
\begin{equation}\label{Riccati-1}
\left \{
\begin{split}
  \begin{aligned}
& \frac{d}{dt}\phi_t+2\rho_t\phi_t+N^{-1}b_t^2\phi_t^2-L=0,\\
  & \phi_T=-R.
      \end{aligned}
  \end{split}
  \right.
\end{equation}
and
\begin{equation}\label{Riccati-2}
\left \{
\begin{split}
  \begin{aligned}
d\varphi_t=-\big(&\big(2\rho_t+2N^{-1}b_t^2\phi_t+N^{-1}b_t^2\EE(\varphi_t)\big)\varphi_t+a_tA_t\EE(\varphi_t)\\
&+N^{-1}b_t^2\EE(A_t\varphi_t)\chi_t+a_t\phi_tA_t\big)dt+\Lambda_tdW_t,\\
d\chi_t=-\big(&\big(2\rho_t+2N^{-1}b_t^2\phi_t+\al_t+a_t\EE(A_t)+N^{-1}b_t^2\EE(A_t\chi_t)\big)\chi_t+a_tA_t\EE(\chi_t)\\
&+(a_2+N^{-1}b_t^2\EE(\chi_t))\varphi_t+a_t\phi_t\big)dt+\Lambda_tdW_t,\\
d\psi_t=-\big(&(\rho_t+N^{-1}b_t^2\phi_t)\psi_t+a_tA_t\EE(\psi_t)+N^{-1}b_t^2\EE(\psi_t)\varphi_t\\
&+(\be c_t\EE(A_t)+N^{-1}b_t^2\EE(A_t\psi_t))\chi_t\big)dt+\Sigma_tdW_t,\\
   \varphi_T&=0\ \ \ \chi_T=0 \ \ \ \psi_T=Rc_0-\lambda.
      \end{aligned}
  \end{split}
  \right.
\end{equation}
Hence,
\begin{equation}\label{feedback form}
u_t=N^{-1}b_t\left(\phi_t X_t+\varphi_t\EE X_t+\chi_t\EE(A_tX_t)+\psi_t\right),
\end{equation}
with  $\phi_t, \varphi_t,\chi_t,\psi_t$, given by Riccati equation (\ref{Riccati-1}) and  the linear BSDE  (\ref{Riccati-2}). 

Note that  $\psi_T=Rc_0-\lambda$  implies that $\EE(\psi_t),\EE(A_t\psi_t) $  are  linear   functions of $\lambda$. We write 
 \[\EE(\psi_t)= \lambda g_t^1+\tilde{g}_t^1 \ \ \ and  \ \ \  \EE(A_t\psi_t)=\lambda g_t^2+\tilde{g}_t^2,\] 
where  $g_t^i$, $\tilde{g}_t^i\; (i=1,2)$  are determined by equations  (\ref{weighted process }) and (\ref{Riccati-2}) explicitly.

To obtain the optimal premium policy, we must identify the value of $\lambda$. To his end,  pluging (\ref{feedback form}) into (\ref{eq0441 }), and  set
\[\xi_t^1=\rho_t+N^{-1}b_t^2\phi_t+N^{-1}b_t^2\EE(\varphi_t),\]
 \[\xi_t^2=a_t+N^{-1}b_t^2\EE(\chi_t),\] 
\[\eta_t^1=N^{-1}b_t^2\EE(A_t\varphi_t),\] 
\[\eta_t^2=\rho_t+a_t\EE(A_t)+\al_t+N^{-1}b_t^2\phi_t+N^{-1}b_t^2\EE(A_t\chi_t),\]
 we obtain  linear ODEs:
\begin{equation}\label{example3-1}
\left \{
\begin{split}
  \begin{aligned}
d\EE(X_t)&=\big(\xi_t^1\EE(X_t)+\xi_t^2\EE(A_tX_t)+N^{-1}b_t^2\EE(\psi_t)\big)dt\\
d\EE(A_tX_t&)=\big(\eta_t^1\EE(X_t)+\eta_t^2\EE(A_tX_t)+N^{-1}b_t^2\EE(A_t\psi_t)+c_t\be_t\EE(A_t)\big)dt\\
\EE X_0&=x_0,\ \ \ \ \EE (A_0X_0)=a_0x_0,
     \end{aligned}
  \end{split}
  \right.
\end{equation} 
Since $\xi_t^i,\eta_t^i,i=1,2$ are continuous on an interval $[0,T]$ and can  be determined by equation (\ref{weighted process }), (\ref{eq0441 }) and  (\ref{Riccati-2}),  quations (\ref{example3-1}) have  a unique explicit solution
\[M_t=M_0\exp\left(\int_0^tF_sds\right)+\int_0^tN_s(\lambda)e^{(\int_0^sF_rdr)}ds:=\lambda(h_t^1,h_t^2)^*+(\tilde {h}_t^1,\tilde{h}_t^2)^*,\]
where $M_t=(\EE(X_t), \EE(A_tX_t))^*,M_0=(x_0, a_0x_0)^* $, $F_t=(\zeta^1_t,\zeta^2_t)$ with $\zeta^1_t=(\xi^1_t,\eta^1_t)^*$, $\zeta^2_t=(\xi^2_t,\eta^2_t)^* $ and 
 \[N_s(\lambda)= \lambda N^{-1}b_t^2(g_t^1,g_t^2)^*+N^{-1}b_t^2(\tilde{g}_t^1,\tilde{g}_t^2+Nb_t^{-2}c_t\be_t\EE(A_t))^*,\] 
 $A^*$ denotes the transpose of the matrix $A$.  Finally  we can  identify $\lambda$   by $c_0=\EE X_T=\lambda h_T^1+\tilde{h}_T^1$. 
 
\begin{remark} The model studied  above is inspired by Wang and Wu \cite{W-W-2022}, where the liability prcess depend on the average cash flow. In our example we consider a more realistic model to allow the  liability prcess to depend on the average of weighted cash flow with a terminal constraint. 
\end{remark}

\section{The proof of the existence and uniqueness}

In this section, we present the proof of Theorem \ref{thm0730a}. To this end, we need the following lemma.

\begin{lemma}
Under (\ref{2})  of Hypothesis \ref{hyp}, $\forall\ p\ge 1$, we have
\[\sup_{t\le T}\EE(A_t^p)<\infty.\]
\end{lemma}
Proof: By It\^o's formula, it is easy to show that
\begin{equation}\label{Lemma-1}
A_t=a\exp\(\int^t_0\tilde{\al}(X_s,\mu_s,u_s)ds+\int^t_0\be(X_s,\mu_s,u_s)dW_s\)\end{equation}
where $\tilde{\al}=\al-\frac12\be^2$. Then,
\begin{eqnarray*}
\EE(A_t^p)&=&a^p\EE\exp\(\int^t_0p\tilde{\al}(X_s,\mu_s,u_s)ds+\int^t_0p\be(X_s,\mu_s,u_s)dW_s\)\\
&\le&a^p\tilde{\EE}\exp\(\int^t_0\(p\tilde{\al}+\frac{p^2}{2}\be^2\)(X_s,\mu_s,u_s)ds\)\\
&\le&K,
\end{eqnarray*}
where $\tilde{\EE}$ is the expectation with respect to an equivalent probability measure given by Girsanov's formula.
\qed

Now, we are ready to present

{\em Proof of Theorem \ref{thm0730a}} :
 We denote $\tilde{A}_t=\ln A_t$. For the existence, we take a Picard sequence with
\[X^{n+1}_t=x+\int^t_0b(X^n_s,\mu^{X^n,A^n}_s,u_s)ds+\int^t_0\si(X^n_s,\mu^{X^n,A^n}_s,u_s)dW_s\]
and
\[\tilde{A}^{n+1}_t=\ln a+\int^t_0\tilde{\al}(X^n_s,\mu^{X^n,A^n}_s,u_s)ds+\int^t_0\be(X^n_s,\mu^{X^n,A^n}_s,u_s)dW_s.\]
 Then,
\begin{eqnarray*}
\EE\sup_{s\le t}|X^{n+1}_s-X^n_s|^2
&\le&K\int^t_0\(\EE|X^n_s-X^{n-1}_s|^2+\rho(\mu^{X^n,A^n}_s,\mu^{X^{n-1},A^{n-1}}_s)^2\)ds.
\end{eqnarray*}
As
\begin{eqnarray*}
\rho(\mu^{X^n,A^n}_s,\mu^{X^{n-1},A^{n-1}}_s)
&=&\sup_{f\in\BB_1}\left|\EE(A^n_sf(X^n_s)-A^{n-1}_sf(X^{n-1}_s))\right|\\
&\le&\EE|A^n_s-A^{n-1}_s|+\EE\(A^{n-1}_s|X^n_s-X^{n-1}_s|\)\\
&\le&\EE|A^n_s-A^{n-1}_s|+K\(\EE|X^n_s-X^{n-1}_s|^2\)^{1/2},
\end{eqnarray*}
we can continue with
\[\EE\sup_{s\le t}|X^{n+1}_s-X^n_s|^2\le K\int^t_0\(\EE|X^n_s-X^{n-1}_s|^2+\(\EE|A^n_s-A^{n-1}_s|\)^2\)ds.\]
Similarly, we have
\[\EE\sup_{s\le t}|\tilde{A}^{n+1}_s-\tilde{A}^n_s|^2\le K\int^t_0\(\EE|X^n_s-X^{n-1}_s|^2+\(\EE|A^n_s-A^{n-1}_s|\)^2\)ds.\]
Then,
\begin{eqnarray*}
\(\EE\sup_{s\le t}|A^{n+1}_s-A^n_s|\)^2
&\le&\(\EE\sup_{s\le t}(A^{n+1}_s+A^n_s)|\tilde{A}^{n+1}_s-\tilde{A}^n|\)^2\\
&\le&K\EE\sup_{s\le t}|\tilde{A}^{n+1}_s-\tilde{A}^n|^2\\
&\le&K\int^t_0\(\EE|X^n_s-X^{n-1}_s|^2+\(\EE|A^n_s-A^{n-1}_s|\)^2\)ds.\end{eqnarray*}
Let
\[f^n(t)\equiv\EE\sup_{s\le t}|X^{n+1}_s-X^n_s|^2+\(\EE\sup_{s\le t}|A^{n+1}_s-A^n_s|\)^2.\]
Then,
\begin{equation}\label{eq0613b}
f^n(t)\le K\int^t_0f^{n-1}(s)ds.\end{equation}
By iterating, we see that
\[f^n(T)\le K'\frac{(KT)^n}{n!}\]
which is summarable. Hence, there exists process $(X,A)$ such that
\[\EE\sup_{s\le T}|X^n_s-X_s|^2+\(\EE\sup_{s\le t}|A^n_s-A_s|\)^2\to 0.\]
It is then easy to show that $(X,A)$ is a solution to (\ref{eq0613a-}).

To prove the uniqueness, we take two solutions and denote the difference of them as $(Y,B)$. Let
\[g(t)\equiv\EE\sup_{s\le t}|Y_s|^2+\(\EE\sup_{s\le t}|B_s|\)^2.\]
Similar to (\ref{eq0613b}), we get
\[g(t)\le K\int^t_0g(s)ds\]
and hence $g=0$. This yields the uniqueness.
\qed

\section{The proof of the stochastic maximum principle}

Suppose that $u_t$ is an optimal control and $v_t$ is such that $v_t+u_t\in\mathcal{U}$. Then, $u^\ep_t\equiv u_t+\ep v_t\in\mathcal{U}$. Hence,
$J(u)\le J(u^\ep)$ for all $\ep\in[0,1]$.  Let $(X^\ep,A^\ep)$ be the state with control $u^\ep$. Then,
\[\left\{\begin{array}{ccl}
dX^\ep_t&=&b(X^\ep_t,\mu^\ep_t,u^\ep_t)dt+\si(X_t,\mu^\ep_t,u^\ep_t)dW_t\\
dA^\ep_t&=&A^\ep_t\(\al(X^\ep_t,\mu^\ep_t,u^\ep_t)dt+\be(X^\ep_t,\mu^\ep_t,u^\ep_t)dW_t\)\\
X^\ep_0&=&x,\ A^\ep_0=a.
\end{array}\right.\]
Define $Y^\ep_t=X^\ep_t-X_t$. Then,
\[dY^\ep_t=\(b(X^\ep_t,\mu^\ep_t,u^\ep_t)-b(X_t,\mu_t,u_t)\)dt+\(\si(X_t^\ep,\mu^\ep_t,u^\ep_t)-\si(X_t,\mu_t,u_t)\)dW_t.\]
Let $\th^\ep_t$ be between $(X^\ep_t,\mu^\ep_t,u^\ep_t)$ and $(X_t,\mu_t,u_t)$ such that
\begin{equation}\label{eq0611f}
b(X^\ep_t,\mu^\ep_t,u^\ep_t)-b(X_t,\mu_t,u_t)=b_x(\th^\ep_t)Y^\ep_t+\<\mu^\ep_t-\mu_t,b_\mu(\th^\ep_t,\cdot)\>+\ep b_u(\th^\ep_t)v_t.\end{equation}
Note that
\begin{eqnarray}\label{eq0611g}
|\<\mu^\ep_t-\mu_t,b_\mu(\th^\ep_t,\cdot)\>|
&=&\left|\EE\(A^\ep_tb_\mu(\th^\ep_t,X^\ep_t)-A_tb_\mu(\th^\ep_t,X_t)\)\right|\nonumber\\
&\le&K\EE|A^\ep_t-A_t|+K\EE\(A_t|Y^\ep_t|\).
\end{eqnarray}
 We first give the following lemmas.
\begin{lemma}Under the Hypothesis \ref{hyp} on the coefficients we have, 
\begin{equation}\label{eq0610a}
\EE|X^\ep_t-X_t|^2+\(\EE|A^\ep_t-A_t|\)^2\le K\ep^2.\end{equation}
\end{lemma}
Proof: By (\ref{eq0611f}, \ref{eq0611g}), we obtain
\[|b(X^\ep_t,\mu^\ep_t,u^\ep_t)-b(X_t,\mu_t,u_t)|^2
\le K|Y^\ep_t|^2+K\(\EE|A^\ep_t-A_t|\)^2+K\EE(|Y^\ep_t|^2)+K\ep^2.\]
The same estimate holds for $\si$. Then,
\begin{eqnarray}\label{eq0611h}
\EE(|Y^\ep_t|^2)&\le&K\EE\(\int^t_0|b(X^\ep_s,\mu^\ep_s,u^\ep_s)-b(X_s,\mu_s,u_s)|ds\)^2\nonumber\\
&&+K\EE\int^t_0|\si(X^\ep_s,\mu^\ep_s,u^\ep_s)-\si(X_s,\mu_s,u_s)|^2ds\nonumber\\
&\le&K\int^t_0\(\EE(|Y^\ep_s|^2)+\(\EE|A^\ep_s-A_s|\)^2+\ep^2\)ds.
\end{eqnarray}

On the other hand, we denote $\tilde{A}_t=\ln A_t$. Then,
\begin{eqnarray}\label{eq0611i}
\(\EE|A^\ep_t-A_t|\)^2
&=&\(\EE\left|{\ze^\ep_t}\(\int^t_0\(\tilde{\al}(X^\ep_s,\mu^\ep_s,u^\ep_s)-\tilde{\al}(X_s,\mu_s,u_s)\)ds\right.\right.\right.\nonumber\\
&&\qquad\left.\left.\left.+\int^t_0\(\be(X^\ep_s,\mu^\ep_s,u^\ep_s)-\be(X_s,\mu_s,u_s)\)dW_s\)\right|\)^2\nonumber\\
&\le&K\int^t_0\(\EE(|Y^\ep_s|^2)+\(\EE|A^\ep_s-A_s|\)^2+\ep^2\)ds,
\end{eqnarray}
where $\ze^\ep_t$ is between $A^\ep_t$ and $A_t$. 
Denote the left hand side  of (\ref{eq0610a}) by $f(t)$. Adding inequalities (\ref{eq0611h}, \ref{eq0611i}), we obtain
\[f(t)\le K\int^t_0 f(s)ds+K\ep^2.\]
The desired conclusion follows from Gronwall's inequality.
\qed

Next, we hope to prove that
\[\lim_{\ep\to 0}\ep^{-1}(X^\ep_t-X_t)=Y_t\mbox{ and }\lim_{\ep\to 0}\ep^{-1}(A^\ep_t-A_t)=B_t,\]
where $Y_t$ and $B_t$ are two processes. We first taking derivative formally to guess the form of $Y_t$ and $B_t$. In fact, we will show that
\begin{equation}
\left \{
\begin{split}
  \begin{aligned}
dY_t&=\Big(b_x(X_t,\mu_t,u_t)Y_t+b_u(X_t,\mu_t,u_t)v_t\\
& \ \ \ +\EE'\(B_t'b_{\mu}(X_t,\mu_t,u_t; X_t')+A_t'b_{\mu,1}(X_t,\mu_t,u_t; X_t')Y_t'\)\Big)dt\\
& \ \ \ +\Big(\si_x(X_t,\mu_t,u_t)Y_t+\si_u(X_t,\mu_t,u_t)v_t\\
& \ \ \ +\EE'\(B_t'\si_{\mu}(X_t,\mu_t,u_t; X_t')+A_t'\si_{\mu,1}(X_t,\mu_t,u_t; X_t')Y_t'\)\Big)dW_t,\\
Y_0&=0,
 \end{aligned}
  \end{split}
  \right.
\end{equation}
and 

\begin{equation}
\left \{
\begin{split}
  \begin{aligned}
d\tilde{B}_t&=\Big(\tilde{\al}_x(X_t,\mu_t,u_t)Y_t+\tilde{\al}_u(X_t,\mu_t,u_t)v_t\\
  & \ \ \ +\EE'\(B_t'\tilde{\al}_{\mu}(X_t,\mu_t,u_t; X_t')+A_t'\tilde{\al}_{\mu,1}(X_t,\mu_t,u_t; X_t')Y_t'\)\Big)dt\\
& \ \ \ +\Big(\be_x(X_t,\mu_t,u_t)Y_t+\be_u(X_t,\mu_t,u_t)v_t\\
& \ \ \ +\EE'\(B_t'\be_{\mu}(X_t,\mu_t,u_t; X_t')+A_t'\be_{\mu,1}(X_t,\mu_t,u_t; X_t')Y_t'\)\Big)dW_t,\\
B_0&=0,
\end{aligned}
  \end{split}
  \right.
\end{equation}

where $B_t=A_t\tilde{B}_t$.

\begin{lemma}Under  Hypothesis \ref{hyp} we have,
\[\lim_{\ep\to 0}\ep^{-1}(X^\ep_t-X_t)=Y_t\mbox{ and }\lim_{\ep\to 0}\ep^{-1}(A^\ep_t-A_t)=B_t.\]
\end{lemma}
Proof: Let
\[Z^\ep_t=\ep^{-1}(X^\ep_t-X_t)-Y_t\mbox{ and }C^\ep_t=\ep^{-1}(A^\ep_t-A_t)-B_t.\]
Then
\[dZ^\ep_t=b^\ep(t)dt+\si^\ep(t)dW_t,\]
where
\begin{eqnarray*}
b^\ep(t)&=&\ep^{-1}\big(b(X^\ep_t,\mu^\ep_t,u^\ep_t)-b(X_t,\mu_t,u_t)\big)
-\big(b_x(X_t,\mu_t,u_t)Y_t+b_u(X_t,\mu_t,u_t)v_t\big)\\
&&-\EE'\big(B_t'b_{\mu}(X_t,\mu_t,u_t; X_t')+A_t'b_{\mu,1}(X_t,\mu_t,u_t; X_t')Y_t'\big),
\end{eqnarray*}
and $\si^\ep(t)$ is given similarly.

Note that
\begin{eqnarray*}
b^\ep(t)&=&b_x(\th^\ep_t)\ep^{-1}Y^\ep_t-b_x(X_t,\mu_t,u_t)Y_t+\ep^{-1}\<\mu^\ep_t-\mu_t,b_\mu(\th^\ep_t,\cdot)\>\\
&&-\EE'\big(B_t'b_{\mu}(X_t,\mu_t,u_t; X_t')+A_t'b_{\mu,1}(X_t,\mu_t,u_t; X_t')Y_t'\big)\\
&&+ b_u(\th^\ep_t)v_t-b_u(X_t,\mu_t,u_t)v_t\\
&=&b_x(\th^\ep_t)Z^\ep_t+\de b_x(t)Y_t+\de b_u(t)v_t\\
&&+\ep^{-1}\EE'\big(A'^\ep_tb_\mu(\th^\ep_t; X'^\ep_t)-A'_tb_\mu(\th^\ep_t; X'_t)\big)\\
&&-\EE'\big(B_t'b_{\mu}(X_t,\mu_t,u_t; X_t')+A_t'b_{\mu,1}(X_t,\mu_t,u_t; X_t')Y_t'\big)\\
&=&b_x(\th^\ep_t)Z^\ep_t+\de b_x(t)Y_t+\de b_u(t)v_t+\EE'\big(C'^\ep_tb_\mu(\th^\ep_t; X'^\ep_t)+Z'^\ep_tA'_tb_{\mu,1}(\th^\ep_t; \ze'^\ep_t)\big)\\
&&+\EE'\big(A_t'Y'_t\(b_{\mu,1}(\th^\ep_t; \ze'^\ep_t)-b_{\mu,1}(\th_t; X'_t)\)+B'_t\( b_\mu(\th^\ep_t; X'^\ep_t)-b_\mu(\th_t; X'_t)\)\big)
\end{eqnarray*}
where
\[\de b_x(t)=b_x(\th^\ep_t)-b_x(X_t,\mu_t,u_t)\to 0.\]
Similar estimate holds for $\si^\ep(t)$. Then,
\begin{equation}\label{eq0611a}
\EE|Z^\ep_t|^2\le K\int^t_0\(\EE|Z^\ep_s|^2+(\EE|C^\ep_s|)^{2}+\EE\de^\ep_s\)ds,\end{equation}
where $\EE\de^\ep_s\to 0$.

On the other hand, we denote
\[\tilde{C}^\ep_t=\ep^{-1}(\tilde{A}^\ep_t-\tilde{A}_t)-\tilde{B}_t.\]
Then
\[C^\ep_t=A_t\tilde{C}^\ep_t+A_t\ep^{-1}\(e^{\ep(\tilde{B}_t+\tilde{C}^\ep_t)}-1-\ep(\tilde{B}_t+\tilde{C}^\ep_t)\)
\equiv A_t\tilde{C}^\ep_t+\de^\ep_t,\]
and $\EE|\de^\ep_t|\to 0$.

Similar to above, we have
\[\EE|\tilde{C}^\ep_t|^2\le K\int^t_0\(\EE|Z^\ep_s|^2+(\EE|C^\ep_s|)^{2}+\EE\de^\ep_s\)ds.\]
Thus,
\begin{eqnarray}\label{eq0611b}
(\EE|C^\ep_t|)^{2}&\le& K\EE|\tilde{C}^\ep_t|^2+(\EE|\de^\ep_t|)^{2}\nonumber\\
&\le&K\int^t_0\(\EE|Z^\ep_s|^2+(\EE|C^\ep_s|)^{2}+\EE\de^\ep_s\)ds+(\EE|\de^\ep_t|)^{2}.\end{eqnarray}
The conclusion follows from (\ref{eq0611a}, \ref{eq0611b}).
\qed

 Denote $\ka_t=(X_t,A_t,\mu_t,u_t)$ and note that as $\ep\to 0$, 
\begin{eqnarray}\label{eq0611c}
0&\le&\ep^{-1}\(J(u^\ep)-J(u)\)\nonumber\\
&\to&\EE\int^T_0\(f_x(\ka_t)Y_t+f_a(\ka_t)B_t+\EE'\(B_t'f_\mu(\ka_t; X_t')+A_t'f_{\mu,1}(\ka_t; X_t')Y_t'\)+f_u(\ka_t)v_t\)dt\nonumber\\
&&+\EE\(\Phi_x(X_T,A_T)Y_T+\Phi_a(X_T,A_T)B_T\).
\end{eqnarray}
Recall “adjoint processes” $(p,q), (P,Q)$ statisties equations (\ref{eq0613adjointa1}), (\ref{eq0613adjointa2}), which we rewrite as  (\ref{eq0613adjointa1}), (\ref{eq0613adjointa2}) as
\[dp_t=g_tdt+q_tdW_t, p_T=-\Phi_a(X_T,A_T),\qquad dP_t=G_tdt+Q_tdW_t, P_T=-\Phi_x(X_T,A_T),\]
where
\begin{equation}\label{eqp}
  \begin{aligned}
  g_t=&-\big\{p_t\al(\th_t)+q_t\be(\th_t)-f_a(\ka_t)\big\}&\\
    &-\mathbb{E}' \big\{(P_t'b_{\mu}+Q_t'\si_{\mu})(\th_t';X_t)+A_t'(p_t'\al_{\mu}+q_t'\be_{\mu})(\th_t';X_t)-f_\mu(\ka_t';X_t)\big\},&
 \end{aligned}
\end{equation}
\begin{equation}\label{eqP}
  \begin{aligned}
   G_t=&-\big\{P_tb_x(\th_t)+Q_t\si_x(\th_t)+A_t(p_t\al_x(\th_t)+q_t\be_x(\th_t))-f_x(\ka_t)\big\}&\\
    &-\mathbb{E}' \big\{A_t\big(P_t'b_{\mu1}(\th_t';X_t)+Q_t'\si_{\mu1}(\th_t';X_t)\big)\big\}&\\
&-\mathbb{E}' \big\{A_t\big(A_t'(p_t'\al_{\mu1}(\th_t';X_t)+q_t'\be_{\mu1}(\th_t';X_t))-f_{\mu1}(\ka_t';X_t)\big)\big\}&
      \end{aligned}
\end{equation}
where  $\th_t=(X_t,\mu_t,u_t)$. Applying It\^o's formula of $p_tB_t$, we get
\begin{eqnarray*}
d(p_tB_t)&=&p_tdB_t+B_tdp_t+d\<B,p\>_t.
\end{eqnarray*}
To continue, we calculate
\begin{eqnarray*}
dB_t&=&A_td\tilde{B}_t+\tilde{B}_tdA_t+d\<A,\tilde{B}\>_t\\
&=&A_t\(\tilde{\al}_x(\th_t)Y_t+\EE'\(B_t'\tilde{\al}_{\mu}(\th_t; X_t')+A_t'\tilde{\al}_{\mu,1}(\th_t; X_t')Y_t'\)+\tilde{\al}_u(\th_t)v_t\)dt\\
&&+A_t\(\be_x(\th_t)Y_t+\EE'\(B_t'\be_{\mu}(\th_t; X_t')+A_t'\be_{\mu,1}(\th_t; X_t')Y_t'\)+\be_u(\th_t)v_t\)dW_t\\
&&+\tilde{B}_tA_t\(\al(\th_t)dt+\be(\th_t)dW_t\)\\
&&+\(\be_x(\th_t)Y_t+\EE'\(B_t'\be_{\mu}(\th_t; X_t')+A_t'\be_{\mu,1}(\th_t; X_t')Y_t'\)+\be_u(\th_t)v_t\)A_t\be(\th_t)dt\\
&=&\(A_t\al_x(\th_t)Y_t+A_t\al_u(\th_t)v_t+B_t\al(\th_t)
\)dt\\
&&+A_t\EE'\(B_t'\al_{\mu}(\th_t; X_t')+A_t'Y_t'\al_{\mu,1}(\th_t; X_t')\)dt\\
&&+\(A_t\(\be_x(\th_t)Y_t+\be_u(\th_t)v_t\)+B_t\be(\th_t)\)dW_t\\
&&+A_t\EE'\(B_t'\be_{\mu}(\th_t; X_t')+A_t'\be_{\mu,1}(\th_t; X_t')Y_t'\)dW_t
\end{eqnarray*}
Thus,
\begin{eqnarray*}
d(p_tB_t)&=&p_t\(A_t\al_x(\th_t)Y_t+A_t\al_u(\th_t)v_t+B_t\al(\th_t)
\)dt\\
&&+p_tA_t\EE'\(B_t'\al_{\mu}(\th_t; X_t')+A_t'Y_t'\al_{\mu,1}(\th_t; X_t')\)dt
+B_tg_tdt\\
&&+q_t\(A_t\(\be_x(\th_t)Y_t+\be_u(\th_t)v_t\)+B_t\be(\th_t)\)dt\\
&&+q_tA_t \EE'\(B_t'\be_{\mu}(\th_t; X_t')+A_t'\be_{\mu,1}(\th_t; X_t')Y_t'\)dt\\
&&+(\cdots)dW_t.
\end{eqnarray*}
Taking integration and expectation, we get
\begin{eqnarray}\label{eq0611d}
&&-\EE(\Phi_a(X_T,A_T)B_T)\\
&=&\EE\int^T_0p_t\(A_t\al_x(\th_t)Y_t+A_t\al_u(\th_t)v_t+B_t\al(\th_t)
\)dt\nonumber\\
&&+\EE\int^T_0\(p_tA_t\EE'\(B_t'\al_{\mu}(\th_t; X_t')+A_t'Y_t'\al_{\mu,1}(\th_t; X_t')\)
+B_tg_t\)dt\nonumber\\
&&+\EE\int^T_0q_t\(A_t\(\be_x(\th_t)Y_t+\be_u(\th_t)v_t\)+B_t\be(\th_t)\)dt\nonumber\\
&&+\EE\int^T_0q_tA_t\EE'\(B_t'\be_{\mu}(\th_t; X_t')+A_t'\be_{\mu,1}(\th_t; X_t')Y_t'\)dt.\nonumber
\end{eqnarray}
Similarly, applying It\^o's formula to $P_tY_t$, we get
\begin{eqnarray*}
d(P_tY_t)&=&P_tdY_t+Y_tdP_t+d\<Y,P\>_t\\
&=&P_t\(b_x(\th_t)Y_t+\EE'\(B_t'b_{\mu}(\th_t; X_t')+A_t'b_{\mu,1}(\th_t; X_t')Y_t'\)+b_u(\th_t)v_t\)dt\\
&&+P_t\(\si_x(\th_t)Y_t+\EE'\(B_t'\si_{\mu}(\th_t; X_t')+A_t'\si_{\mu,1}(\th_t; X_t')Y_t'\)+\si_u(\th_t)v_t\)dW_t\\
&&+Y_tG_tdt+Y_tQ_tdW_t\\
&&+Q_t\(\si_x(\th_t)Y_t+\EE'\(B_t'\si_{\mu}(\th_t; X_t')+A_t'\si_{\mu,1}(\th_t; X_t')Y_t'\)+\si_u(\th_t)v_t\)dt\\
&=&\(P_t\(b_x(\th_t)Y_t+b_u(\th_t)v_t\)+Q_t\(\si_x(\th_t)Y_t+\si_u(\th_t)v_t\)+Y_tG_t\)dt\\
&&+\EE'\(B_t'(P_tb_{\mu}(\th_t; X_t')+Q_t\si_{\mu}(\th_t; X_t'))\)dt\\
&&+\EE'\(A_t'Y_t'(P_tb_{\mu,1}(\th_t; X_t')+Q_t\si_{\mu,1}(\th_t; X_t'))\)dt\\  
&&+(\cdots)dW_t.
\end{eqnarray*}
Taking integration and expectation, we get
\begin{eqnarray}\label{eq0611e}
&&-\EE(\Phi_x(X_T,A_T)Y_T)\\
&=&\EE\int^T_0\(P_t\(b_x(\th_t)Y_t+b_u(\th_t)v_t\)+Q_t\(\si_x(\th_t)Y_t+\si_u(\th_t)v_t\)+Y_tG_t\)dt\nonumber\\
&&+\EE\int^T_0\EE'\big(B_t'\big(P_tb_{\mu}(\th_t; X_t')+Q_t\si_{\mu}(\th_t; X_t')\big)+A_t'Y_t'\big(P_tb_{\mu,1}(\th_t; X_t')+Q_t\si_{\mu,1}(\th_t; X_t')\big)\big)dt.\nonumber
\end{eqnarray}

Combining (\ref{eq0611c}-\ref{eq0611e}), we get
\begin{eqnarray*}
0&\le&\EE\int^T_0\big(f_x(\ka_t)Y_t+f_a(\ka_t)B_t+\EE'\(B_t'f_\mu(\ka_t; X_t')+A_t'f_{\mu,1}(\ka_t; X_t')Y_t'\)+f_u(\ka_t)v_t\big)dt\\
&&+\EE\int^T_0\(-P_t\(b_x(\th_t)Y_t+b_u(\th_t)v_t\)-Q_t\(\si_x(\th_t)Y_t+\si_u(\th_t)v_t\)-Y_tG_t\)dt\\
&&+\EE\int^T_0\EE'\big(-B_t'(P_tb_{\mu}(\th_t; X_t')+Q_t\si_{\mu}(\th_t; X_t'))\big)dt\\
&&+\EE\int^T_0\EE'\big(-A_t'Y_t'(P_tb_{\mu,1}(\th_t; X_t')+Q_t\si_{\mu,1}(\th_t; X_t'))\big)dt\\
&&+\EE\int^T_0-p_t\big(A_t\al_x(\th_t)Y_t+A_t\al_u(\th_t)v_t+B_t\al(\th_t)
\big)dt\\
&&+\EE\int^T_0\big(-p_tA_t\EE'\(B_t'\al_{\mu}(\th_t; X_t')+A_t'Y_t'\al_{\mu,1}(\th_t; X_t')\)
-B_tg_t\big)dt\\
&&+\EE\int^T_0-q_t\big(A_t\(\be_x(\th_t)Y_t+\be_u(\th_t)v_t\)-B_t\be(\th_t)\big)dt\\
&&+\EE\int^T_0-q_tA_t\EE'\big(B_t'\be_{\mu}(\th_t; X_t')+A_t'\be_{\mu,1}(\th_t; X_t')Y_t'\big)dt\\
&\le-&\EE\int^T_0 v_t\left(P_tb_u(\th_t)+Q_t\si_u(\th_t)+A_t(p_t\al_u(\th_t)+q_t\be_u(\th_t))-f_u(\ka_t)\right)dt
\end{eqnarray*}

By the define of Hamiltion $H(X,A,\mu,u,P,Q,p,q)$, we then have
\[0\le -\EE\int^T_0 v_tH_u(X_t,A_t,\mu_t,u_t,P_t,Q_t,p_t,q_t)dt.\]
This implies
\[H_u(X_t,A_t,\mu_t,u_t,P_t,Q_t,p_t,q_t)=0,\]
and hence, finished the proof.

 \end{document}